\documentclass[oneside,english]{amsart}
\usepackage{mathpazo}
\usepackage[T1]{fontenc}
\usepackage[latin9]{inputenc}
\usepackage{geometry}
\usepackage{varioref}
\usepackage{calc}
\usepackage{amsthm}
\usepackage{graphicx}
\usepackage{amssymb}

\providecommand{\tabularnewline}{\\}

\numberwithin{equation}{section} 
\numberwithin{figure}{section} 
\theoremstyle{plain}
\theoremstyle{plain}
\newtheorem{thm}{Theorem}
 \theoremstyle{definition}
  \newtheorem{example}[thm]{Example}
  \theoremstyle{plain}
  \newtheorem*{thm*}{Theorem}

\usepackage[pdftex]{color}
\definecolor{commentcolor}{rgb}{0.6, 0.6, 0.6}

\newcommand\seswarrows{\ensuremath{\textrm{\rlap{$\swarrow$}$\searrow$}}}

\newcommand\closer{\ \scalebox{0.4}{\begin{picture}(0,0)
    \thicklines
    \put(0,10){\line(1,-1){10}}
    \put(0,10){\line(1,1){10}}
    \put(10,0){\line(1,1){10}}
    \put(10,20){\line(1,-1){10}}
    \put(15,0){\line(1,1){10}}
    \put(15,20){\line(1,-1){10}}
  \end{picture}}}

\usepackage{babel}

\begin{document}

\title[Generalizing Dodgson's method]{Generalizing Dodgson's method: a ``double-crossing'' approach to
computing determinants}

\author{}
\begin{abstract}
Dodgson's method of computing determinants was recently revisited
in a paper that appeared in the \emph{College Math Journal}. The method
is attractive, but fails if an interior entry of an intermediate matrix
has the value zero. This paper reviews the structure of Dodgson's
method and introduces a generalization, called a ``double-crossing''
method, that provides a workaround to the failure for many interesting
cases.
\end{abstract}
\maketitle

\section{\label{sec: Introduction}Introduction}

Algebra students learn this simple pattern for the determinant of
a $2\times2$ matrix:\[
\left|\begin{array}{rcr}
a &  & b\\
 & \seswarrows\\
c &  & d\end{array}\right|=ad-bc.\]
A similar pattern exists for $3\times3$ matrices, but to compute
the determinant of larger matrices, the student must learn something
a little more complicated. Most students learn to compute determinants
by expansion of minors, first developed by Laplace. Some students
also learn to compute determinants by triangularizing the matrix.
Both methods are effective, and triangularization is efficient, but
\begin{itemize}
\item hand computations frequently lead to many mistakes;
\item expansion of minors is tedious; and
\item triangularization can turn a matrix of integers into a matrix of fractions.
\end{itemize}
In 1866, the Rev. Charles Lutwidge Dodgson%
\footnote{Dodgson is better known to children of all ages as Lewis Carroll,
author of ``Alice in Wonderland'' and ``Jabberwocky''.%
}, developed a conceptually simple method to compute determinants~\cite{Shutting Up}.
Dodgson's method iterates the familiar $2\times2$ formula; although
it uses division, matrices with integer entries do not turn into matrices
with rational entries.
\begin{example}
\label{exa: basic example of Dodgson's method}Given the matrix\[
A=\left(\begin{array}{rrr}
1 & 0 & 1\\
1 & 3 & 1\\
0 & 1 & 1\end{array}\right),\]
we set $A^{\left(3\right)}=A$. To compute $A^{\left(2\right)}$,
compute the determinants of the $2\times2$ contiguous submatrices
of $A^{\left(3\right)}$:\[
A^{\left(2\right)}=\left(\begin{array}{cc}
\left|\begin{array}{rr}
1 & 0\\
1 & 3\end{array}\right| & \left|\begin{array}{rr}
0 & 1\\
3 & 1\end{array}\right|\\
\\\left|\begin{array}{rr}
1 & 3\\
0 & 1\end{array}\right| & \left|\begin{array}{rr}
3 & 1\\
1 & 1\end{array}\right|\end{array}\right)=\left(\begin{array}{rr}
3 & -3\\
1 & 2\end{array}\right).\]
To compute $A^{\left(1\right)}$, repeat the computation for $A^{\left(2\right)}$,
but divide by the center entry of $A^{\left(3\right)}$:\[
A^{\left(1\right)}=\left(\frac{\left|\begin{array}{rr}
3 & -3\\
1 & 2\end{array}\right|}{3}\right)=\left(\frac{9}{3}\right)=\left(3\right).\]
Now go back and compute $\det A$ using your favorite method (pattern,
expansion of cofactors, triangularization, etc.). You will find that
$\det A=3$.\closer
\end{example}
Dodgson's method is \emph{quick} and \emph{conceptually simple}. In
general, we can describe the method in the following way:
\begin{itemize}
\item Let $A^{\left(n\right)}$ be the $n\times n$ matrix given.
\item For each $k=n-1,n-2,\ldots,1$:

\begin{itemize}
\item Let $B^{\left(k\right)}$ be the $k\times k$ matrix of determinants
of contiguous $2\times2$ submatrices of $A^{\left(k+1\right)}$.
\item If $k=n-1$, let $A^{\left(k\right)}=B^{\left(k\right)}$.
\item If $k\leq n-2$, let $A^{\left(k\right)}$ be the $k\times k$ matrix
whose $\left(i,j\right)$-th element is the $\left(i,j\right)$-th
element of $B^{\left(k\right)}$ divided by the $\left(i,j\right)$-th
element of the interior of $A^{\left(k+2\right)}$. (The \textbf{interior}
of an $n\times n$ matrix $M$ is the $\left(n-2\right)\times\left(n-2\right)$
submatrix whose $\left(i,j\right)$-th element is the $\left(i+1,j+1\right)$-th
element of $M$).
\end{itemize}
\item The singleton element of $A^{\left(1\right)}$ is the determinant
of $A$.
\end{itemize}
Dodgon's method is an example of a \textbf{condensation} method to
compute determinants; each iterate $A^{\left(k\right)}$ is a condensation
of the previous iterate $A^{\left(k+1\right)}$. Other condensation
methods appear in~\cite{Aitken}. If Dodgson's method terminates
successfully, it computes the determinant of an $n\times n$ matrix
using\[
4\left[\left(n-1\right)^{2}+\left(n-2\right)^{2}+\cdots+1^{2}\right]\approx4\cdot n\cdot n^{2}\]
multiplications, subtractions, and divisions, or $O\left(n^{3}\right)$
operations in $\mathbb{Q}$. This is not bad, especially considering
that all the divisions are exact, so no fractions are introduced;
Bareiss' algorithm, a better-known fraction-free method of computing
determinants, also performs $O\left(n^{3}\right)$ operations in $\mathbb{Q}$.~\cite{Modern Computer Algebra,Yap2000}

There's the rub, though: division presents Dodgson's method with a
huge drawback. (It also presents an obstacle with the usage of Bareiss'
algorithm.) What's so bad about division?
\begin{example}
\label{exa: swap rows of good basic ex}Swap the first two rows of
the matrix of Example~\ref{exa: basic example of Dodgson's method}
to obtain\[
M=\left(\begin{array}{rrr}
1 & 3 & 1\\
1 & 0 & 1\\
0 & 1 & 1\end{array}\right),\]
We know from the properties of determinants that $\det M=-\det A=-3$.

What happens when we use Dodgson's method? First we compute\[
M^{\left(2\right)}=\left(\begin{array}{rr}
-3 & 3\\
1 & -1\end{array}\right).\]
To compute $M^{\left(1\right)}$ we divide the determinant of $M^{\left(2\right)}$
by the interior element of $M^{\left(3\right)}$. \emph{But the interior
element of $M^{\left(3\right)}$ is zero!}\closer
\end{example}
In general, Dodgson's method fails to compute the determinant of a
matrix $A$ whenever a zero appears in the interior of $A^{\left(k\right)}$
for any $k\geq3$. This can happen even if no zeroes appear in the
interior of $A$.

A workaround discussed in \cite{Shutting Up} swaps rows of the \emph{original}
matrix in such a way that zeroes are moved out of the interior. For
example, if you swap the top two rows of $M$ in Example~\ref{exa: swap rows of good basic ex},
you return to $A$ of example~\ref{exa: basic example of Dodgson's method},
for which Dodgson's method worked fine. However, there are two drawbacks
to this workaround. \emph{First,} swapping rows may well introduce
other zeroes into the matrix, and it isn't easy to predict this from
the outset. \emph{Second,} swapping rows simply won't work for some
matrices.
\begin{example}
\label{exa: unswappable}No combination of row or column swaps will
allow Dodgson's method to compute the determinant of\[
N=\left(\begin{array}{rrrr}
1 & 0 & 3 & 0\\
0 & -1 & 0 & 1\\
1 & 1 & 2 & 0\\
0 & 2 & 0 & 1\end{array}\right),\]
because there will \emph{always} be a zero in the interior of $N$.\closer
\end{example}
Will a different workaround of Dodgson's method work for $N$? \emph{Yes!}
We describe such a method in Section~\ref{sec: The fix}. The reader
will see quickly why we call it a ``double-crossing'' method.

\section{\label{sec: The fix}A ``double-crossing'' method}

Before describing the method, we illustrate it using the matrix from
Example~\vref{exa: swap rows of good basic ex}.
\begin{example}
\label{exa: fix, initial 0}Recall \begin{equation}
M^{\left(3\right)}=\left(\begin{array}{rrr}
1 & 3 & 1\\
1 & 0 & 1\\
0 & 1 & 1\end{array}\right)\qquad\mathrm{and}\qquad M^{\left(2\right)}=\left(\begin{array}{rr}
-3 & 3\\
1 & -1\end{array}\right).\label{eq: M2 normal}\end{equation}
The zero in the interior of $M^{\left(3\right)}$ causes Dodgson's
method to fail when computing $M^{\left(1\right)}$.

Above that zero is a non-zero element, $M_{1,2}^{\left(3\right)}=3$.
We will divide by this element instead, but this requires us to re-compute
$M^{\left(2\right)}$ in a slightly different manner. Cross out the
first row and second column of $M^{\left(3\right)}$ (the ones containing
3). We are left with the $2\times2$ complementary matrix\[
M^{*}=\left(\begin{array}{rr}
1 & 1\\
0 & 1\end{array}\right).\]
Recall that we compute the \emph{minor of an element} of a matrix
by (again) crossing out the row and column containing that element,
then taking the determinant of the remaining submatrix. Consider the
matrix $M'$ of \emph{minors of elements} of $M^{\left(3\right)}$
that correspond to the elements of $M^{*}$; that is,\[
M'=\left(\begin{array}{rr}
\left|\begin{array}{rr}
3 & 1\\
1 & 1\end{array}\right| & \left|\begin{array}{rr}
1 & 3\\
0 & 1\end{array}\right|\\
\\\left|\begin{array}{rr}
3 & 1\\
0 & 1\end{array}\right| & \left|\begin{array}{rr}
1 & 3\\
1 & 0\end{array}\right|\end{array}\right).\]
Put\begin{equation}
M^{\left(2\right)}=M'=\left(\begin{array}{rr}
2 & 1\\
3 & -3\end{array}\right);\label{eq: M2 double-crossed}\end{equation}
conclude by dividing $\det M^{\left(2\right)}$ by the non-zero element
of $M^{\left(3\right)}$ that we identified earlier: \[
M^{\left(1\right)}=\left(\frac{\left|M^{\left(2\right)}\right|}{3}\right)=\left(\frac{-6-3}{3}\right)=\left(-\frac{9}{3}\right)=\left(-3\right).\]
As noted in example~\ref{exa: swap rows of good basic ex}, $\det M=-3$.

Notice, by the way, that after a row and column swap all but one of
the values in $M^{\left(2\right)}$ in line~\eqref{eq: M2 double-crossed}
are the values of $M^{\left(2\right)}$ in line~\eqref{eq: M2 normal}.
We will say more about this later.\closer
\end{example}
The generalization we have presented of Dodgson's method preserves
its ``spirit'', inasmuch as we computed all determinants by condensing
$2\times2$ contiguous submatrices. The example shows why we call
it the ``double-crossing'' method:
\begin{itemize}
\item we found a non-zero element adjacent to the zero element;
\item we crossed out its row and column, obtaining the complementary matrix
$M^{*}$;
\item for each element in $M^{*}$, we compute its minor by

\begin{itemize}
\item (again) crossing out the row and column of its location in $M^{\left(3\right)}$,
and
\item using the determinant of the remaining matrix to compute the element
of $M^{\left(2\right)}$.
\end{itemize}
\end{itemize}
\noindent When a zero appears in an intermediate matrix $A^{\left(k\right)}$
(where $n>k\geq3$), we cross out rows and columns that make up a
\emph{submatrix} of $A^{\left(n\right)}$. Let's look at a matrix
where a zero does not appear in the interior of $A$, but does appear
in the interior of an intermediate matrix.
\begin{example}
\label{exa: fix, intermediate 0 (5x5)}Let\[
A=\left(\begin{array}{rrrrr}
1 & 0 & 1 & 0 & 1\\
0 & 5 & 3 & 1 & 0\\
1 & 3 & 2 & 1 & 1\\
0 & 1 & 1 & 1 & 0\\
2 & 0 & 2 & 0 & 1\end{array}\right).\]
Using Dodgson's method, we compute $A^{\left(5\right)}$, $A^{\left(4\right)}$,
and then\[
A^{\left(3\right)}=\left(\begin{array}{rrr}
-4 & -2 & 2\\
-2 & 0 & -2\\
4 & -4 & -1\end{array}\right).\]
We encounter a zero in the interior! Above it is a non-zero element,
$-2$. Notice that it lies in the top row and central column of $A^{\left(3\right)}$.

To compute $A^{\left(2\right)}$, do the following:
\begin{itemize}
\item cross out the top three rows and central three columns of $A^{\left(5\right)}$
(the \emph{original} matrix);
\item this gives us a $2\times2$ complementary matrix $M^{*}$;
\item for each element in $M^{*}$,

\begin{itemize}
\item cross out the row and column of its location in $A^{\left(5\right)}$,
\item compute the resulting minor, and
\item put it into the corresponding location in $A^{\left(2\right)}$.
\end{itemize}
\end{itemize}
\noindent Following these instructions, we have\[
M^{*}=\left(\begin{array}{rr}
0 & 0\\
2 & 1\end{array}\right)\qquad\mathrm{and}\qquad A^{\left(2\right)}=\left(\begin{array}{rr}
\left|\begin{array}{rrrr}
0 & 1 & 0 & 1\\
5 & 3 & 1 & 0\\
3 & 2 & 1 & 1\\
0 & 2 & 0 & 1\end{array}\right| & \left|\begin{array}{rrrr}
1 & 0 & 1 & 0\\
0 & 5 & 3 & 1\\
1 & 3 & 2 & 1\\
2 & 0 & 2 & 0\end{array}\right|\\
\\\left|\begin{array}{rrrr}
0 & 1 & 0 & 1\\
5 & 3 & 1 & 0\\
3 & 2 & 1 & 1\\
1 & 1 & 1 & 0\end{array}\right| & \left|\begin{array}{rrrr}
1 & 0 & 1 & 0\\
0 & 5 & 3 & 1\\
1 & 3 & 2 & 1\\
0 & 1 & 1 & 1\end{array}\right|\end{array}\right).\]
We compute the determinants for $A^{\left(2\right)}$ using Dodgson's
method---and again, the bottom two values turn out to have values
that you would obtain in the two top rows by ordinary condensation
of $A^{\left(3\right)}$, while the others are different:\[
A^{\left(2\right)}=\left(\begin{array}{rr}
2 & 0\\
4 & -4\end{array}\right).\]
We now compute $A^{\left(1\right)}$ by dividing the determinant of
$A^{\left(2\right)}$ by the non-zero entry of $A^{\left(3\right)}$
that we identified above: -2. We obtain\[
A^{\left(1\right)}=\left(\frac{\left|\begin{array}{rr}
2 & 0\\
4 & -4\end{array}\right|}{-2}\right)=\left(\frac{-8}{-2}\right)=\left(4\right),\]
and in fact $\det A=4$.\closer
\end{example}
So far we've been using a special case of the double-crossing method.
Theorem~\ref{thm: double-cross, special case} describes this special
case precisely.
\begin{thm}
[Double-crossing method, special case]\label{thm: double-cross, special case}Let
$A$ be an $n\times n$ matrix. Suppose that we try to evaluate $\det A$
using Dodgson's method, but we encounter a zero in the interior of
$A^{\left(k\right)}$, say in row $i$ and column $j$ of $A^{\left(k\right)}$.
If the element $\alpha$ in row $i-1$ and column $j$ of $A^{\left(k\right)}$
is non-zero, then we compute $A^{\left(k-1\right)}$and $A^{\left(k-2\right)}$
as usual, with the following exception for the element in row $i-1$
and column $j-1$ of $A^{\left(k-2\right)}$:
\begin{itemize}
\item let $\ell=n-k$;
\item identify the $\left(\ell+3\right)\times\left(\ell+3\right)$ submatrix
$\mathcal{A}$ whose upper left corner is the element in row $i-1$
and column $j-1$ of $A^{\left(n\right)}$;
\item identify the $2\times2$ complementary minor $M^{*}$ by crossing
out the $\left(\ell+1\right)\times\left(\ell+1\right)$ submatrix
of \textup{$\mathcal{A}$} whose upper left corner is the element
in row~1 and column~2 of $\mathcal{A}$;
\item compute the matrix $M'$ of determinants of minors of elements of
$M^{*}$ in $\mathcal{A}$;
\item compute the element in row $i-1$ and column $j-1$ of $A^{\left(k-2\right)}$
by dividing the determinant of $M'$ by $\alpha$. \closer
\end{itemize}
\end{thm}
A proof of correctness appears in Section~\ref{sec: the fix works!},
after we explain why Dodgson's original method works correctly.

One can generalize Theorem~\ref{thm: double-cross, special case}
so that the non-zero element appears immediately above, below, left,
right, or catty-corner to the zero; that is, the non-zero element
is \emph{adjacent} to the zero: see Theorem~\vref{thm: double-crossing method}.
If the zero appears in a $3\times3$ block of zeroes, then the double-crossing
method will not repair Dodgson's method, although additional strategies
may be possible.

Before proceeding to the next section, we encourage the reader to
go back and examine how Theorem~\ref{thm: double-cross, special case}
describes what we did in the examples of this section. We conclude
by applying the double-crossing method to compute the determinant
$N$ from Example~\ref{exa: unswappable}.
\begin{example}
\label{exa: double-cross on unswappable}Recall from Example~\ref{exa: unswappable}\[
N=\left(\begin{array}{rrrr}
1 & 0 & 3 & 0\\
0 & -1 & 0 & 1\\
1 & 1 & 2 & 0\\
0 & 2 & 0 & 1\end{array}\right).\]
Put $N^{\left(4\right)}=N$; we have one zero element in the interior,
at $\left(i,j\right)=\left(2,3\right)$. As described in Theorem~\ref{thm: double-cross, special case},
this corresponds to element $\left(1,2\right)$ of $N^{\left(2\right)}$.
The other three elements of $N^{\left(2\right)}$ can be computed
as usual:\[
N^{\left(3\right)}=\left(\begin{array}{rrr}
-1 & 3 & 3\\
1 & -2 & -2\\
2 & -4 & 2\end{array}\right)\quad\longrightarrow\quad N^{\left(2\right)}=\left(\begin{array}{rr}
1 & ?\\
0 & -6\end{array}\right).\]
We use the double-crossing method to compute the other element of
$N^{\left(2\right)}$. The interior zero appears in the original matrix,
so $\ell=0$. We choose the $3\times3$ submatrix whose upper left
corner is the element in row $i-1=1$, column $j-1=2$ of $N^{\left(4\right)}$:\[
\mathcal{A}=\left(\begin{array}{rrr}
0 & 3 & 0\\
-1 & 0 & 1\\
1 & 2 & 0\end{array}\right).\]
Cross out the $1\times1$ submatrix in row~1, column~2 of $\mathcal{A}$
and identify the complementary matrix\[
M^{*}=\left(\begin{array}{rr}
-1 & 1\\
1 & 0\end{array}\right).\]
Compute the corresponding matrix of determinants of minors in $\mathcal{A}$\[
M'=\left(\begin{array}{cc}
\left|\begin{array}{rr}
3 & 0\\
2 & 0\end{array}\right| & \left|\begin{array}{rr}
0 & 3\\
1 & 2\end{array}\right|\\
\\\left|\begin{array}{rr}
3 & 0\\
0 & 1\end{array}\right| & \left|\begin{array}{rr}
0 & 3\\
-1 & 0\end{array}\right|\end{array}\right)=\left(\begin{array}{rr}
0 & -3\\
3 & 3\end{array}\right).\]
The determinant of this matrix is 9; after dividing by the nonzero
$N_{1,3}^{\left(4\right)}=3$ we have\[
N^{\left(2\right)}=\left(\begin{array}{rr}
1 & 3\\
0 & -6\end{array}\right).\]

We can now conclude by computing\[
N^{\left(1\right)}=\left(\frac{\left|\begin{array}{rr}
1 & 3\\
0 & -6\end{array}\right|}{-2}\right)=\left(3\right).\]
In fact, the determinant of $N$ is 3.\closer
\end{example}
The most computationally intensive part of the double-crossing method
is that of computing the matrix of determinants of minors, but many
of these \emph{have already been computed}. It likely is not clear
to the reader at present how to identify these, but as we explore
the mechanics of Dodgson's method and the double-crossing method,
we will see that we can predict exactly which parts of the matrix
of cofactors need recomputing, and which can be copied from previous
work.

In any case, it is time to consider why Dodgson's method works.

\section{\label{sec: Jacobi's Theorem}Why does Dodgson's method work?}

Just as Bareiss' algorithm relies on a well-known theorem of Sylvester,
Dodgson's method relies on the following theorem of Jacobi.~\cite{Shutting Up}
\begin{thm*}
[Jacobi's Theorem]Let
\begin{itemize}
\item $\mathcal{A}$ be an $n\times n$ matrix;
\item $M$ an $m\times m$ minor of $\mathcal{A}$, where $m<n$, chosen
from rows $i_{1},i_{2},\ldots,i_{m}$ and columns $j_{1},j_{2},\ldots,j_{m}$,
\item $M'$ the corresponding $m\times m$ minor of $\mathcal{A}'$, the
matrix of cofactors of $\mathcal{A}$, and
\item $M^{*}$ the $\left(n-m\right)\times\left(n-m\right)$
minor of $\mathcal{A}$ complementary to $M$.
\end{itemize}
Then\[
\det M'=\left(\det\mathcal{A}\right)^{m-1}\cdot\det M^{*}\cdot\left(-1\right)^{\sum_{\ell=1}^{m}i_{\ell}+j_{\ell}}.\mbox{\closer}\]

\end{thm*}
Section~\ref{sec: the fix works!} gives a proof of the new method.

We adopt the following notation. If $A\in\mathbb{R}^{m\times n}$
denotes a matrix, then $A'$ denotes its matrix of cofactors. Both
Dodgson's method and the double-crossing method make use of submatrices
of a matrix, which for a given matrix $A$ we denote in the following
fashion:

\begin{center}
\begin{tabular}{cl}
$A_{i\ldots j,k\ldots\ell}$ & rows $i,i+1,\ldots,j$ and columns $k,k+1,\ldots,\ell$ of $A$;\tabularnewline
$A_{i,k\ldots\ell}$ & row $i$ and columns $k,k+1,\ldots,\ell$ of $A$;\tabularnewline
$A_{i\ldots j,k}$ & rows $i,i+1,\ldots,j$ and column $k$ of $A$.\tabularnewline
\end{tabular}
\par\end{center}
\begin{example}
Recall from Example~\ref{exa: swap rows of good basic ex} \[
A=\left(\begin{array}{rrr}
1 & 3 & 1\\
1 & 0 & 1\\
0 & 1 & 1\end{array}\right).\]
Its matrix of cofactors is\[
A'=\left(\begin{array}{rcrcr}
\left|\begin{array}{rr}
0 & 1\\
1 & 1\end{array}\right| &  & -\left|\begin{array}{rr}
1 & 1\\
0 & 1\end{array}\right| &  & \left|\begin{array}{rr}
1 & 0\\
0 & 1\end{array}\right|\\
\\-\left|\begin{array}{rr}
3 & 1\\
1 & 1\end{array}\right| &  & \left|\begin{array}{rr}
1 & 1\\
0 & 1\end{array}\right| &  & -\left|\begin{array}{rr}
1 & 3\\
0 & 1\end{array}\right|\\
\\\left|\begin{array}{rr}
3 & 1\\
0 & 1\end{array}\right| &  & -\left|\begin{array}{rr}
1 & 1\\
1 & 1\end{array}\right| &  & \left|\begin{array}{rr}
1 & 3\\
1 & 0\end{array}\right|\end{array}\right)\]
while\[
A_{1\ldots2,1\ldots2}=\left(\begin{array}{rr}
1 & 3\\
1 & 0\end{array}\right)\qquad\mathrm{and}\qquad A_{3,2\ldots3}=\left(\begin{array}{rr}
1 & 1\end{array}\right).\mbox{\closer}\]

\end{example}
To illustrate the relationship between Dodgson's method and Jacobi's
Theorem, consider a generic $4\times4$ matrix.%
\begin{figure}
\caption{\label{fig:Diagram-for-4x4}Diagram illustrating how Jacobi's Theorem
is applied to Dodgson's method, given the $4\times4$ matrix of Example~\ref{exa: illustration of relationship}.}

\centering{}%
\framebox{\begin{minipage}[t]{0.95\columnwidth}%
\begin{center}
\includegraphics[scale=0.2]{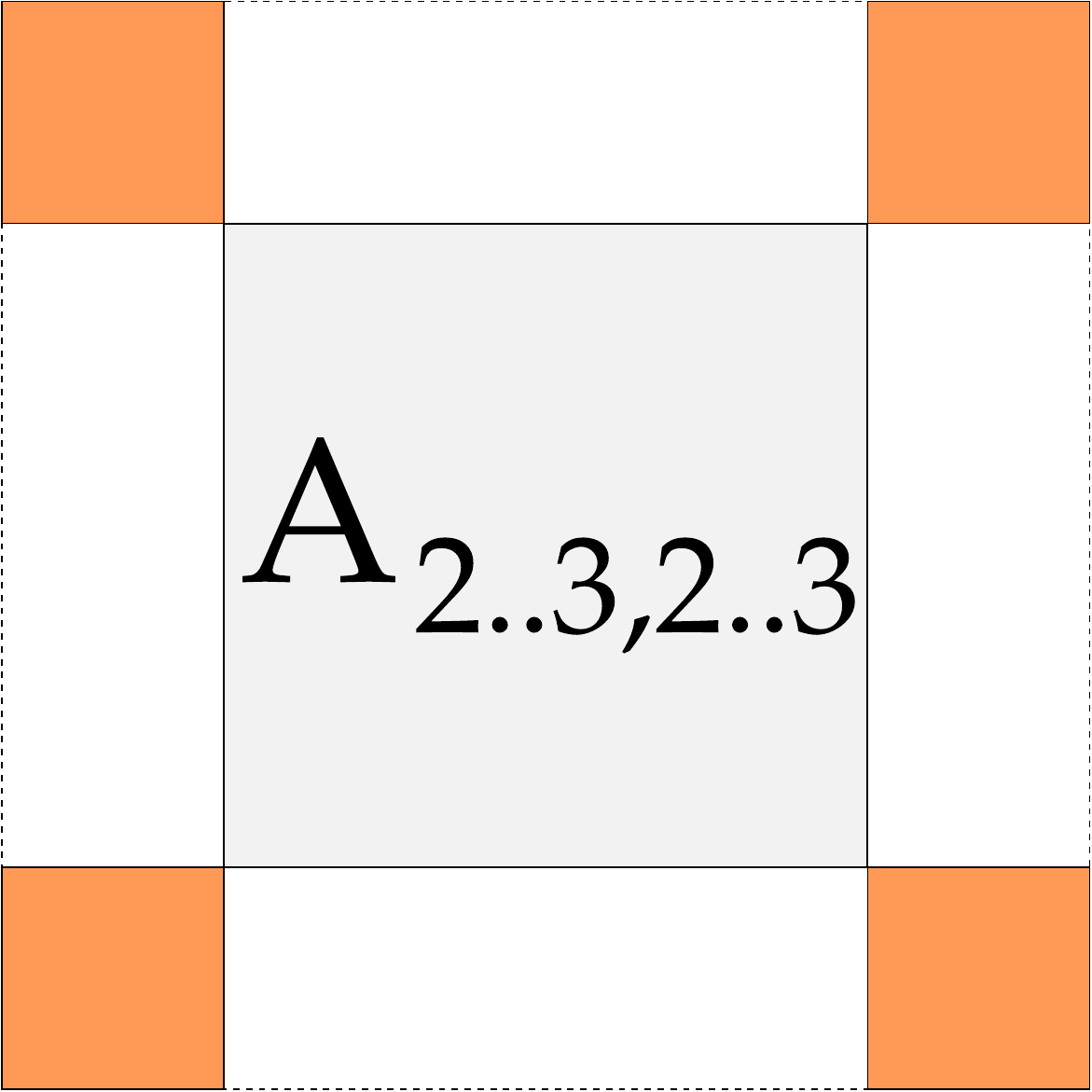}
\par\end{center}

Consider the generic $4\times4$ matrix from Example~\ref{exa: illustration of relationship}.
The final condensation from Dodgson's method gave us\[
A^{\left(1\right)}=\left(\frac{|A_{1\ldots3,1\ldots3}||A_{2\ldots4,2\ldots4}|-|A_{2\ldots4,1\ldots3}||A_{1\ldots3,2\ldots4}|}{|A_{2\ldots3,2\ldots3}|}\right).\]
For Jacobi's Theorem, select the $2\times2$ minor of the corners
of $A$\[
M=\left(\begin{array}{rr}
A_{1,1} & A_{1,4}\\
A_{4,1} & A_{4,4}\end{array}\right).\]
Its complementary minor is $M^{*}=\left(A_{2\ldots3,2\ldots3}\right)$.
The matrix of cofactors is\[
M'=\left(\begin{array}{rr}
A_{1,1}' & -A_{1,4}'\\
-A_{4,1}' & A_{4,4}'\end{array}\right)\]
where $A_{1,1}'=\left|A_{2\ldots4,2\ldots4}\right|$, $A_{4,4}'=\left|A_{1\ldots3,1\ldots3}\right|$,
$A_{1,4}'=\left|A_{2\ldots4,1\ldots3}\right|$, and $A_{4,1}'=\left|A_{1\ldots3,2\ldots4}\right|$.
By Jacobi's Theorem,\begin{align*}
\det M' & =\left(\det A\right)^{2-1}\cdot\det M^{*}\cdot\left(-1\right)^{1+1+4+4}\\
\frac{\det M'}{\det M^{*}} & =\det A\\
\frac{\left|A_{1\ldots3,1\ldots3}\right|\left|A_{2\ldots4,2\ldots4}\right|-\left|A_{2\ldots4,1\ldots3}\right|\left|A_{1\ldots3,2\ldots4}\right|}{\left|A_{2\ldots3,2\ldots3}\right|} & =\det A.\end{align*}
This is precisely the final step in Dodgson's method as applied to
$A$. Notice that the negatives in $M'$ cancel when computing the
determinant.%
\end{minipage}}
\end{figure}

\begin{example}
\label{exa: illustration of relationship}Let $A$ be a generic $4\times4$
matrix; we show how Dodgson's method applies Jacobi's Theorem. The
first two condensations by Dodgson's method produce\begin{align}
A^{\left(3\right)} & =\left(\begin{array}{ccc}
|A_{1...2,1...2}| & |A_{1...2,2...3}| & |A_{1...2,3...4}|\\
\\|A_{2...3,1...2}| & |A_{2...3,2...3}| & |A_{2...3,3...4}|\\
\\|A_{3...4,1...2}| & |A_{3...4,2...3}| & |A_{3...4,3...4}|\end{array}\right)\nonumber \\
 & \mathrm{and}\nonumber \\
A^{\left(2\right)} & =\left(\begin{array}{cc}
\frac{\left|\begin{array}{cc}
|A_{1...2,1...2}| & |A_{1...2,2...3}|\\
|A_{2...3,1...2}| & |A_{2...3,2...3}|\end{array}\right|}{a_{2,2}} & \frac{\left|\begin{array}{cc}
|A_{1...2,2...3}| & |A_{1...2,3...4}|\\
|A_{2...3,2...3}| & |A_{2...3,3...4}|\end{array}\right|}{a_{2,3}}\\
\frac{\left|\begin{array}{cc}
|A_{2...3,1...2}| & |A_{2...3,2...3}|\\
|A_{3...4,1...2}| & |A_{3...4,2...3}|\end{array}\right|}{a_{3,2}} & \frac{\left|\begin{array}{cc}
|A_{2...3,2...3}| & |A_{2...3,3...4}|\\
|A_{3...4,2...3}| & |A_{3...4,3...4}|\end{array}\right|}{a_{3,3}}\end{array}\right).\label{eq: A(2) in example of Dodgsons relation to Jacobi}\end{align}

To see how $A^{\left(2\right)}$ corresponds to Jacobi's Theorem,
consider the upper left $3\times3$ submatrix of $A^{\left(4\right)}$\[
\mathcal{A}=A_{1\ldots3,1\ldots3}=\left(\begin{array}{ccc}
a_{1,1} & a_{1,2} & a_{1,3}\\
a_{2,1} & a_{2,2} & a_{2,3}\\
a_{3,1} & a_{3,2} & a_{3,3}\end{array}\right).\]
Cross out row~2 and column~2 of $\mathcal{A}$, obtaining\[
M=\left(\begin{array}{rr}
a_{1,1} & a_{1,3}\\
a_{3,1} & a_{3,3}\end{array}\right);\]
its complement in $\mathcal{A}$ is $M^{*}=\left(a_{2,2}\right)$.
The $2\times2$ submatrix of $\mathcal{A}'$ corresponding to the
cofactors of $M$ in $\mathcal{A}$ is\[
M'=\left(\begin{array}{rr}
\left|\begin{array}{rr}
a_{2,2} & a_{2,3}\\
a_{3,2} & a_{3,3}\end{array}\right| & \left|\begin{array}{rr}
a_{2,1} & a_{2,2}\\
a_{3,1} & a_{3,2}\end{array}\right|\\
\\\left|\begin{array}{rr}
a_{1,2} & a_{1,3}\\
a_{2,2} & a_{2,3}\end{array}\right| & \left|\begin{array}{rr}
a_{1,1} & a_{1,2}\\
a_{2,1} & a_{2,2}\end{array}\right|\end{array}\right)=\left(\begin{array}{rr}
\left|A_{2\ldots3,2\ldots3}\right| & \left|A_{2\ldots3,1\ldots2}\right|\\
\left|A_{1\ldots2,2\ldots3}\right| & \left|A_{1\ldots2,1\ldots2}\right|\end{array}\right).\]
Use Jacobi's Theorem and a column and row swap to obtain\begin{align*}
\det M' & =\left(\det\mathcal{A}\right)^{2-1}\cdot\det M^{*}\cdot\left(-1\right)^{1+1+3+3}\\
\frac{\det M'}{\det M^{*}} & =\det\mathcal{A}\\
\frac{\left|\begin{array}{rr}
\left|A_{1\ldots2,1\ldots2}\right| & \left|A_{1\ldots2,2\ldots3}\right|\\
\left|A_{2\ldots3,1\ldots2}\right| & \left|A_{2\ldots3,2\ldots3}\right|\end{array}\right|}{a_{2,2}} & =\det\mathcal{A}.\end{align*}
Since $\mathcal{A}=A_{1\ldots3,1\ldots3}$, we have computed the determinant
of the upper $3\times3$ submatrix of $A$. This is equivalent to
the element in the upper left corner of $A^{\left(2\right)}$ in \ref{eq: A(2) in example of Dodgsons relation to Jacobi};
the negatives from the cofactors in Jacobi's method cancel each other
out. Likewise,
\begin{itemize}
\item the upper right corner of $A^{\left(2\right)}$ has the value $\det A_{1\ldots3,2\ldots4}$;
\item the lower left corner has the value $\det A_{2\ldots4,1\ldots3}$;
and
\item the lower right corner has the value $\det A_{2\ldots4,2\ldots4}$.
\end{itemize}
The next (and final) condensation in Dodgson's method is

\[
A^{\left(1\right)}=\left(\frac{|A^{\left(2\right)}|}{A_{2,2}^{\left(3\right)}}\right)=\left(\frac{|A_{1\ldots3,1\ldots3}||A_{2\ldots4,2\ldots4}|-|A_{2\ldots4,1\ldots3}||A_{1\ldots3,2\ldots4}|}{|A_{2\ldots3,2\ldots3}|}\right).\]
Apply Jacobi's Theorem \emph{on $A$}, using for $M^{*}$ the $2\times2$
interior submatrix of $A$, and we have

\[
|A|=\frac{|A_{1\ldots3,1\ldots3}||A_{2\ldots4,2\ldots4}|-|A_{2\ldots4,1\ldots3}||A_{1\ldots3,2\ldots4}|}{|A_{2\ldots3,2\ldots3}|}.\]
This is precisely the singleton element of $A^{\left(1\right)}$!
See Figure~\ref{fig:Diagram-for-4x4}.\closer\end{example}
\begin{thm}[Dodgson's Condensation Theorem]
Let $A$ be an $n\times n$ matrix. After
$k$ successful condensations, Dodgson's method produces the matrix\[
A^{\left(n-k\right)}=\left(\begin{array}{cccc}
\left|A_{1\ldots k+1,1\ldots k+1}\right| & \left|A_{1\ldots k+1,2\ldots k+2}\right| & \cdots & \left|A_{1\ldots k+1,n-k\ldots n}\right|\\
\left|A_{2\ldots k+2,1\ldots k+1}\right| & \left|A_{2\ldots k+2,2\ldots k+2}\right| & \cdots & \left|A_{2\ldots k+2,n-k\ldots n}\right|\\
\vdots & \vdots & \ddots & \vdots\\
\left|A_{n-k\ldots n,1\ldots k+1}\right| & \left|A_{n-k\ldots n,2\ldots k+2}\right| & \cdots & \left|A_{n-k\ldots n,n-k\ldots n}\right|\end{array}\right)\]
whose entries are the determinants of the $\left(k+1\right)\times\left(k+1\right)$
submatrices of $A$.\closer
\end{thm}
By ``successful'' iterations we mean that one never encounters division
by zero.%
\begin{figure}
\caption{\label{fig:Diagram-for-deltaxdelta}Diagram illustrating how Jacobi's
Theorem is applied to Dodgson's method, given a general $n\times n$
matrix.}

\centering{}%
\framebox{\begin{minipage}[t]{0.95\columnwidth}%
\begin{center}
\includegraphics[scale=0.3]{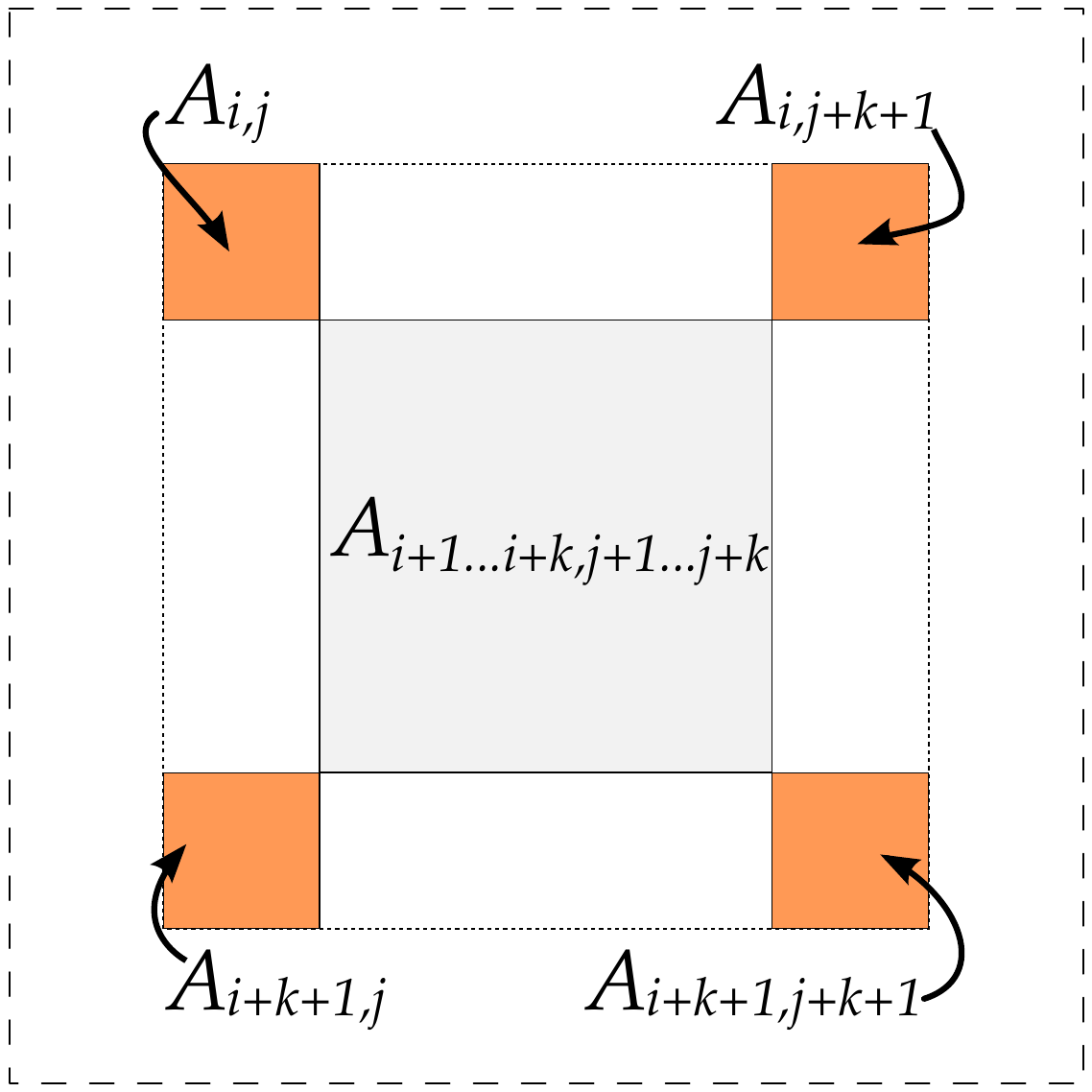}
\par\end{center}

Consider a $\left(k+1\right)\times\left(k+1\right)$ submatrix of
$A$, $\mathcal{A}=A_{i\ldots i+k+1,j\ldots j+k+1}$. Select the $2\times2$
minor \[
M=\left(\begin{array}{cc}
A_{i,j} & A_{i,j+k+1}\\
A_{i+k+1,j} & A_{i+k+1,j+k+1}\end{array}\right).\]
Its complementary minor is $M^{*}=\left(A_{i+1\ldots i+k,j+1\ldots j+k}\right)$.
The matrix of cofactors of $M$ in $\mathcal{A}$ is\[
M'=\left(\begin{array}{cc}
A_{i,j}' & \pm A_{i,j+k+1}'\\
\pm A_{i+k+1,j}' & A_{i+k+1,j+k+1}'\end{array}\right)\]
where\begin{align*}
A_{i,j}' & =\left|A_{i+1\ldots i+k+1,j+1\ldots j+k+1}\right|\\
A_{i,j+k+1}' & =\left|A_{i+1\ldots i+k+1,j\ldots j+k}\right|\\
A_{i+k+1,j}' & =\left|A_{i\ldots i+k,j+1\ldots j+k+1}\right|\\
A_{i+k+1,j+k+1}' & =\left|A_{i\ldots i+k,j\ldots j+k}\right|.\end{align*}
By Jacobi's Theorem,\[
\det M'=\left(\det\mathcal{A}\right)^{2-1}\cdot\det M^{*}\cdot\left(-1\right)^{i+\left(j+k+1\right)+j+\left(i+k+1\right)}\]
or\begin{equation}
\det\mathcal{M}=\frac{\left|A_{i,j}'\right|\cdot\left|A_{i+k+1,j+k+1}'\right|-\left|A_{i,j+k+1}'\right|\cdot\left|A_{i+k+1,j}'\right|}{\left|A_{i+1\ldots i+k,j+1\ldots j+k}\right|},\label{eq: general relationship in fig}\end{equation}
so long as the denominator is not zero.%
\end{minipage}}
\end{figure}

\begin{proof}
We proceed by induction on $k$.

\emph{Inductive Base:} When $k=1$, the theorem is trivial: one condensation
gives\[
A^{\left(n-1\right)}=\left(\begin{array}{cccc}
\left|A_{1\ldots2,1\ldots2}\right| & \left|A_{1\ldots2,2\ldots3}\right| & \cdots & \left|A_{1\ldots2,n-1\ldots n}\right|\\
\left|A_{2\ldots3,1\ldots2}\right| & \left|A_{2\ldots i+2,2\ldots i+2}\right| & \cdots & \left|A_{2\ldots3,n-1\ldots n}\right|\\
\vdots & \vdots & \ddots & \vdots\\
\left|A_{n-1\ldots n,1\ldots2}\right| & \left|A_{n-i\ldots n,2\ldots i+2}\right| & \cdots & \left|A_{n-1\ldots n,n-1\ldots n}\right|\end{array}\right).\]

\emph{Inductive Hypothesis:} Fix $k$. Assume that for all $\ell=1,\ldots,k$,
the $\ell$th condensation gives us $A^{\left(n-\ell\right)}$ where
for all $1\leq i,j\leq n$\[
A_{i,j}^{\left(n-\ell\right)}=\left|A_{i\ldots i+\left(n-\ell\right),j\ldots j+\left(n-\ell\right)}\right|.\]

\emph{Inductive Step:} Let $i,j\in\left\{ 1,\ldots,n-k\right\} $.
The next condensation in Dodgson's method gives us\[
A_{i,j}^{\left(n-(k+1)\right)}=\frac{A_{i,j}^{\left(n-k\right)}A_{i+1,j+1}^{\left(n-k\right)}-A_{i+1,j}^{\left(n-k\right)}A_{i,j+1}^{\left(n-k\right)}}{A_{i+1,j+1}^{\left(n-(k-1)\right)}}.\]
From the inductive hypothesis, we can substitute\begin{align*}
A_{i,j}^{\left(n-(k+1)\right)} & =\frac{\left(\begin{array}{l}
\left|A_{i\ldots i+k,j\ldots j+k}\right|\left|A_{i+1\ldots i+k+1,j+1\ldots j+k+1}\right|\\
\qquad-\left|A_{i+1\ldots i+k+1,j\ldots j+k}\right|\left|A_{i\ldots i+k,j+1\ldots j+k+1}\right|\end{array}\right)}{\left|A_{i+1\ldots i+k,j+1\ldots j+k}\right|}.\end{align*}
Apply Jacobi's Theorem with
\begin{itemize}
\item $\mathcal{A}=A_{i\ldots i+k+1,j\ldots j+k+1}$,
\item $M$ the $2\times2$ minor made up of the corners of $\mathcal{A}$,
\item $M'$ the corresponding $2\times2$ minor of $\mathcal{A}'$, and
\item $M^{*}$ the complementary $k\times k$ minor of $\mathcal{A}$
\end{itemize}
\noindent to see that\[
A_{i,j}^{\left(n-(k+1)\right)}=\left|A_{i\ldots i+k+1,j\ldots j+k+1}\right|.\]
Notice that here $i_1=j_1=1$ and $i_2=j_2=k+2$. See Figure~\ref{fig:Diagram-for-deltaxdelta}.
\end{proof}

\section{\label{sec: the fix works!}Why does the ``double-crossing'' method
work?}

The new workaround is based on Dodgson's Condensatioin Theroem.
The goal of the $k$th condensation is to compute each determinant
$\left|A_{i\ldots i+k,j\ldots j+k}\right|$. Dodgson's method fails
when the corresponding denominator of \eqref{eq: general relationship in fig}
is zero.

However, the fraction of \eqref{eq: general relationship in fig}
is not the only way to apply Jacobi's Theorem. As long as \emph{some}
$\left(k-2\right)\times\left(k-2\right)$ minor of $A_{i\ldots i+k,j\ldots j+k}$
has non-zero determinant, we can still recover, reusing most of the
computations already performed in previous steps of Dodgson's method,
and calculating only a few new minors, again using the same approach
as Dodgson's method. For example, the proof of Theorem~\ref{thm: double-cross, special case}
can be summarized by applying Jacobi's Theorem with
\begin{itemize}
\item $\mathcal{A}$ of Theorem~\ref{thm: double-cross, special case}
standing in for $\mathcal{A}$ of Jacobi's Theorem;
\item the $2\times2$ submatrix $M^{*}$ of Theorem~\ref{thm: double-cross, special case}
standing in for $M^{*}$ of Jacobi's Theorem; and 
\item the matrix of determinants of minors of $M'$ of Theorem~\ref{thm: double-cross, special case}
standing in for $M'$ of Jacobi's Theorem.
\end{itemize}
One difference does require investigation: the negatives in the matrix
of cofactors are now in different places! With Dodgson's Method, we
choose the corners for $M$, which has two important consequences
in Jacobi's Theorem:
\begin{itemize}
\item $i_1=j_1=1$ and $i_2=j_2=k$, so that $\left(-1\right)^{\sum_{\ell=1}^mi_\ell+j_\ell}=\left(-1\right)^{(1+1)+(4+4)}=1$;
and
\item the negatives that appear in the $2\times2$ matrix of cofactors appear
off the main diagonal, cancelling each other out.
\end{itemize}
In the special case of the double-crossing method, the matter is a
little more complicated. We have chosen the top-middle minor for $M$,
so:
\begin{itemize}
\item $\left(-1\right)^{\sum_{\ell=1}^{m}i_{\ell}+j_{\ell}}=\left(-1\right)^{[(k-1)+1]+(k+k)}=\left(-1\right)^k$;
and
\item the determinant of the $2\times2$ matrix of cofactors is\[
(-1)^{[(k-1)+1]+(k+k)}D_{11}D_{22}-(-1)^{[(k-1)+k]+(k+1)}D_{12}D_{21}=(-1)^k\cdot D
\]
where $D_{ij}$ is the minor in the cofactor $M_{ij}'$,
and $D=D_{11}D_{22}-D_{12}D_{21}$.
Thus we can disregard the signs and consider only the determinants of the minors.
\end{itemize}
We use this to revisit Example~\ref{exa: unswappable}, which would
not work with Dodgson's method even after swapping rows or columns.
\begin{example}
\label{exa:illustration of workaround}Recall from Example~\ref{exa: unswappable}\[
N=\left(\begin{array}{rrrr}
1 & 0 & 3 & 0\\
0 & -1 & 0 & 1\\
1 & 1 & 2 & 0\\
0 & 2 & 0 & 1\end{array}\right).\]
\begin{figure}
\caption{\label{fig:Diagram-for-workaround}Diagram for workaround for Example
\ref{exa:illustration of workaround}}

\centering{}%
\framebox{\begin{minipage}[t]{0.95\columnwidth}%
\begin{center}
\includegraphics[scale=0.1]{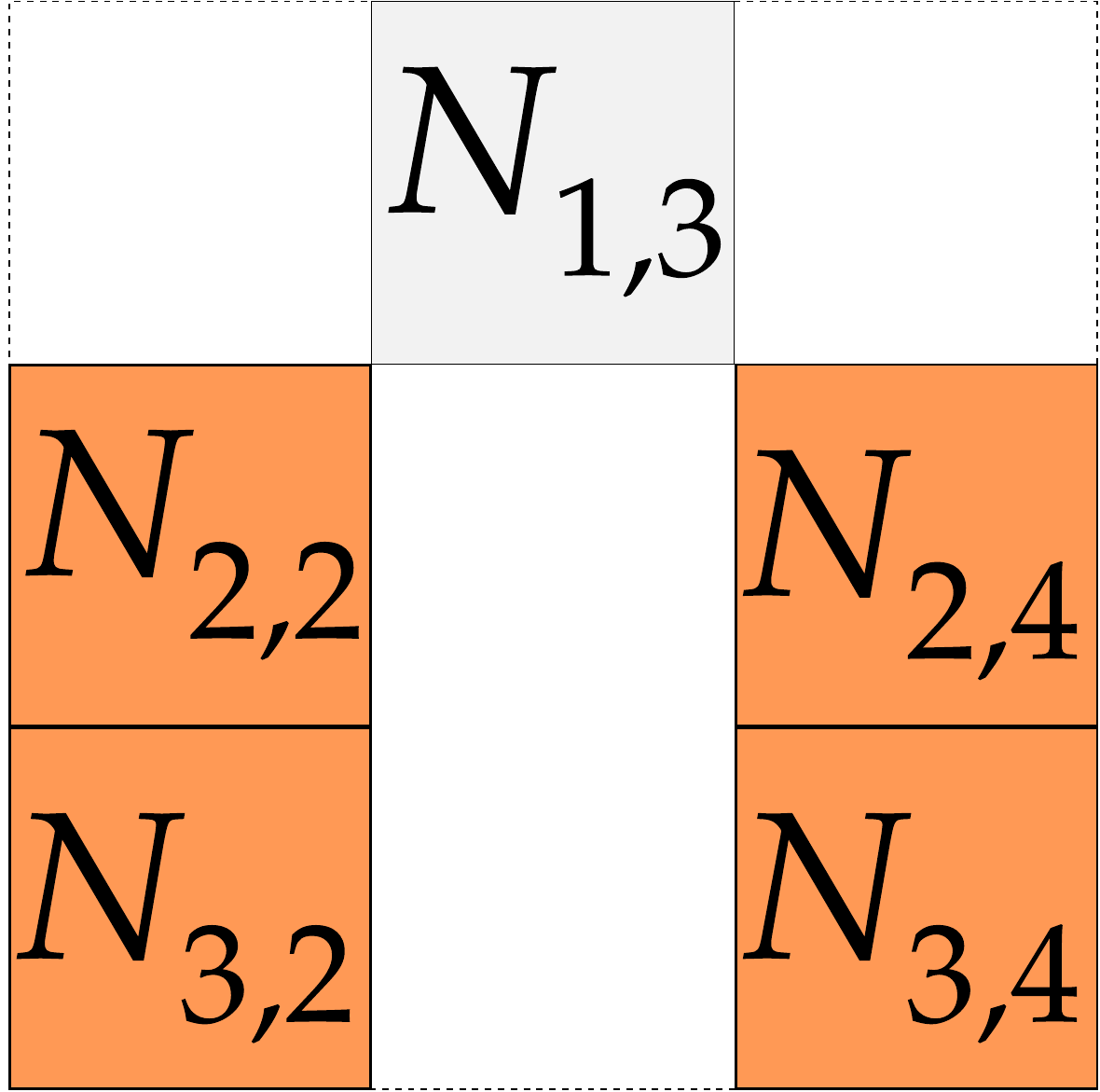}
\par\end{center}

The problem in computing $N_{1,2}^{\left(2\right)}$ was the zero
in position $N_{2,3}$. Since $N_{1,2}^{\left(2\right)}=\det N_{2\ldots,4,2\ldots4}$
choose instead the minor\[
M=\left(\begin{array}{cc}
N_{2,2} & N_{2,4}\\
N_{3,2} & N_{3,4}\end{array}\right)\]
whose complementary minor is $M^{*}=N_{1,3}=3$.
Notice that we have $i_1=2$, $j_1=1$, and $i_2=j_2=3$.
In this case the
matrix of cofactors is\begin{align*}
M' & =\left(\begin{array}{rr}
\left(-1\right)^{2+1}\left|\begin{array}{rr}
N_{1,3} & N_{1,4}\\
N_{3,3} & N_{3,4}\end{array}\right| & \left(-1\right)^{2+3}\left|\begin{array}{rr}
N_{1,2} & N_{1,3}\\
N_{3,2} & N_{3,3}\end{array}\right|\\
\\\left(-1\right)^{3+1}\left|\begin{array}{rr}
N_{1,3} & N_{1,4}\\
N_{2,3} & N_{2,4}\end{array}\right| & \left(-1\right)^{3+3}\left|\begin{array}{rr}
N_{1,2} & N_{1,3}\\
N_{2,2} & N_{2,3}\end{array}\right|\end{array}\right)\\
 & =\left(\begin{array}{rr}
-\left|\begin{array}{rr}
3 & 0\\
2 & 0\end{array}\right| & -\left|\begin{array}{rr}
0 & 3\\
1 & 2\end{array}\right|\\
\\\left|\begin{array}{rr}
3 & 0\\
0 & 1\end{array}\right| & \left|\begin{array}{rr}
0 & 3\\
-1 & 0\end{array}\right|\end{array}\right)=\left(\begin{array}{rr}
0 & 3\\
3 & 3\end{array}\right).\end{align*}
By Jacobi's Theorem,\[
\det N_{2\ldots,4,2\ldots4}=\frac{\det M'}{\det M^{*}\cdot\left(-1\right)^{(2+1)+(3+3)}}=\frac{-9}{-3}=3.\]
Notice the correspondence between $M'$ and the matrix $N^{\left(2\right)}$
obtained in Example~\ref{exa: double-cross on unswappable} using
the double-crossing method: multiplying the top row of $M'$ by~$-1$
gives us $N^{\left(2\right)}$. Dividing $\det M'$ by~$-3=3\cdot\left(-1\right)^{9}$
cancels out the negative introduced into the determinant.%
\end{minipage}}
\end{figure}
Dodgson's method would have us compute $N_{1,2}^{\left(2\right)}$
by choosing for Jacobi's Theorem\[
M=\left(\begin{array}{rr}
N_{1,2} & N_{1,4}\\
N_{3,2} & N_{3,4}\end{array}\right)\qquad\Longrightarrow\qquad M^{*}=\left(N_{2,3}\right),\]
whence\[
N_{1,2}^{\left(2\right)}=\det N_{1\ldots3,2\ldots4}=\frac{\det M'}{\det M^{*}}\]
but $\det M^{*}=N_{2,3}=0$.

The double-crossing method allows us to choose instead\[
M=\left(\begin{array}{rr}
N_{2,2} & N_{2,4}\\
N_{3,2} & N_{3,4}\end{array}\right)\qquad\Longrightarrow\qquad M^{*}=\left(N_{1,3}\right).\]
We have $i_1=2$, $j_1=1$, and $i_2=j_2=3$ (these values are relative to $\mathcal A$, not $N$).
Since $N_{1,3}=1$ we have\begin{align*}
N_{1,2}^{\left(2\right)} & =\det N_{1\ldots3,2\ldots4}\\
 & =\frac{\det M'}{\det M^{*}\cdot\left(-1\right)^{(2+1)+(3+3)}}\\
 & =\frac{M_{1,1}'\cdot M_{2,2}'- M_{1,2}'\cdot M_{2,1}'}{-N_{1,3}}\\
 & =\frac{\left(-1\right)^{2+1}\left|\begin{array}{rr}
3 & 0\\
2 & 0\end{array}\right|\cdot\left(-1\right)^{3+3}\left|\begin{array}{rr}
0 & 3\\
-1 & 0\end{array}\right|-\left(-1\right)^{2+3}\left|\begin{array}{rr}
0 & 3\\
1 & 2\end{array}\right|\cdot\left(-1\right)^{3+1}\left|\begin{array}{rr}
3 & 0\\
0 & 1\end{array}\right|}{-3}\\
 & =\frac{-9}{-3}=3.\end{align*}
See Figure~\ref{fig:Diagram-for-workaround}.\closer
\end{example}
We conclude by showing how Jacobi's Theorem likewise justifies the
``double-crossing'' method in Example~\ref{exa: fix, intermediate 0 (5x5)}.
\begin{example}
Recall from Example~\ref{exa: fix, intermediate 0 (5x5)}\[
A=\left(\begin{array}{rrrrr}
1 & 0 & 1 & 0 & 1\\
0 & 5 & 3 & 1 & 0\\
1 & 3 & 2 & 1 & 1\\
0 & 1 & 1 & 1 & 0\\
2 & 0 & 2 & 0 & 1\end{array}\right).\]
When computing $A^{\left(1\right)}$, Dodgson's method wants us to
cross out the middle three rows and the middle three columns of $A$,
obtaining the minor\[
M=\left(\begin{array}{rr}
A_{1,1} & A_{1,5}\\
A_{5,1} & A_{5,5}\end{array}\right)=\left(\begin{array}{rr}
1 & 1\\
2 & 1\end{array}\right)\]
whose complementary minor is\[
M^{*}=\left(\begin{array}{rrr}
A_{2,2} & A_{2,3} & A_{2,4}\\
A_{3,2} & A_{3,3} & A_{3,4}\\
A_{4,2} & A_{4,3} & A_{4,4}\end{array}\right)=\left(\begin{array}{rrr}
5 & 3 & 1\\
3 & 2 & 1\\
1 & 1 & 1\end{array}\right).\]
Unfortunately, the determinant of this latter matrix is zero, which
corresponds to the zero in the interior of\[
A^{\left(3\right)}=\left(\begin{array}{rrr}
-4 & -2 & 2\\
-2 & 0 & -2\\
4 & -4 & 3\end{array}\right).\]

The double-crossing method suggests instead to cross out the top three
rows and middle three columns of $A$, obtaining the minor\[
M=\left(\begin{array}{rr}
A_{4,1} & A_{4,5}\\
A_{5,1} & A_{5,5}\end{array}\right)=\left(\begin{array}{rr}
0 & 0\\
2 & 1\end{array}\right)\]
whose complementary minor is\[
M^{*}=\left(\begin{array}{rrr}
A_{1,2} & A_{1,3} & A_{1,4}\\
A_{2,2} & A_{2,3} & A_{2,4}\\
A_{3,2} & A_{3,3} & A_{3,4}\end{array}\right)=\left(\begin{array}{rrr}
0 & 1 & 0\\
5 & 3 & 1\\
3 & 2 & 1\end{array}\right).\]
The determinant of $M^{*}$ is $A_{1,2}^{\left(3\right)}=-2$; we
will multiply it by $\left(-1\right)^{1+2+2+3+3+4}=-1$. To compute
$A_{1,2}^{\left(2\right)}$ we need $\det M'$, which requires us
to compute the matrix of cofactors of $M$ in $A^{\left(5\right)}$:\[
M'=\left(\begin{array}{rr}
-\left|\begin{array}{rrrr}
0 & 1 & 0 & 1\\
5 & 3 & 1 & 0\\
3 & 2 & 1 & 1\\
0 & 2 & 0 & 1\end{array}\right| & -\left|\begin{array}{rrrr}
1 & 0 & 1 & 0\\
0 & 5 & 3 & 1\\
1 & 3 & 2 & 1\\
2 & 0 & 2 & 0\end{array}\right|\\
\\\left|\begin{array}{rrrr}
0 & 1 & 0 & 1\\
5 & 3 & 1 & 0\\
3 & 2 & 1 & 1\\
1 & 1 & 1 & 0\end{array}\right| & \left|\begin{array}{rrrr}
1 & 0 & 1 & 0\\
0 & 5 & 3 & 1\\
1 & 3 & 2 & 1\\
0 & 1 & 1 & 1\end{array}\right|\end{array}\right).\]
The top row of negatives introduces a negative into $\det M'$, which
cancels with the $-1$ that is multiplied to $\det M^{*}$, cancelling
each other out. We can thus consider the determinants of minors, rather
than the matrix of cofactors.\closer
\end{example}
Of course, it is possible that both an element of the interior is
zero \emph{and} the element above it is zero. The following theorem
provides the promised generalization of Theorem~\ref{thm: double-cross, special case};
its proof uses Jacobi's method in a manner similar to the proof of
Theorem~\ref{thm: double-cross, special case}. Notice that if $r=s=0$,
then Theorem~\ref{thm: double-crossing method} specializes to Dodgson's
method.
\begin{thm}
\label{thm: double-crossing method}Let $A$ be an $n\times n$ matrix.
Suppose that we try to evaluate $\left|A\right|$ using Dodgson's
method, but we encounter a zero in the interior of $A^{\left(k\right)}$,
say in row $i$ and column $j$ of $A^{\left(k\right)}$. Let $r,s\in\left\{ -1,0,1\right\} $.
If the element $\alpha$ in row $i+r$ and column $j+s$ of $A^{\left(k\right)}$
is non-zero, then we compute $A^{\left(k-2\right)}$ as usual, with
the following exception for the element in row $i-1$ and column $j-1$:
\begin{itemize}
\item let $\ell=n-k$;
\item identify the $\left(\ell+3\right)\times\left(\ell+3\right)$ submatrix
$\mathcal{A}$ whose upper left corner is the element in row $i-1$
and column $j-1$ of $A^{\left(n\right)}$;
\item identify the complementary minor $M^{*}$ by crossing out the $\left(\ell+1\right)\times\left(\ell+1\right)$
submatrix $M$ of $\mathcal{A}$ whose upper left corner is the element
in row~$r+2$ and column~$s+2$ of $\mathcal{A}$;
\item compute the matrix $M'$ of determinants of minors of $M^{*}$ in
$\mathcal{A}$;
\item compute the element in row $i-1$ and column $j-1$ of $A^{\left(k+2\right)}$
by dividing the determinant of $M'$ by $\alpha$. \closer
\end{itemize}
\end{thm}

\section{Concluding remarks}

The double-crossing method succeeds as long as Dodgson's method succeeds,
and with the same number of integer operations, since Dodgson's method
is a special case of the double-crossing method. In the worst case
scenario where Dodgson's method generates a number of interior zeroes,
the double-crossing method requires the computation of at least two
$\left(k-1\right)\times\left(k-1\right)$ determinants for every zero
entry of $A^{\left(k\right)}$. However, this changes the number of
integer operations only by a constant, so the double-crossing method
requires only $O\left(n^{3}\right)$ integer operations. In addition,
the double-crossing method preserves the general simplicity and spirit
of Dodgson's method.

The double-crossing method is not guaranteed to succeed; if an intermediate
matrix contains a $3\times3$ block of zeroes, then the double-crossing
method also fails. Sparse matrices provide an excellent example where
the double-crossing method is an abject failure; consider the identity
matrix of order six or higher.

In those cases where the double-crossing method fails, one can still
preserve the computations that work, and adapt a hybrid with another
method. Theorem~\ref{thm: double-crossing method} tells us the precise
minor $\mathcal{A}$ whose determinant we need; we can compute the
determinant of this minor using another method, substitute its value
into row $i-1$ and column $j-1$ of $A^{\left(k+2\right)}$, and
proceed.


\begin{thebibliography}{4}
\bibitem[1]{Aitken}Alexander Aitken. \emph{Determinants and Matrices.}
Interscience Publishers, 1951.

\bibitem[2]{Shutting Up}Adrian Rice and Eve Torrence. ``Shutting
up like a telescope'': Lewis Carroll's ``Curious Condensation Method
for Evaluating Determinants, \emph{The College Mathematics Journal,}
\textbf{38} (2007) 85--95.

\bibitem[3]{Modern Computer Algebra}Joachim von zur Gathen and J\"urgen
Gerhard. \emph{Modern Computer Algebra.} Cambridge University Press,
1999.

\bibitem[4]{Yap2000}Chee Keng Yap. \emph{Fundamental Problems of
Algorithmic Algebra.} Oxford University Press, 2000.
\end{thebibliography}
\end{document}